\documentclass[10pt]{amsart}
\usepackage{amsmath}
\usepackage{amsfonts}
\usepackage{amsthm}

\theoremstyle{plain}
\newtheorem{thm}{Theorem}[section]
\newtheorem{lem}[thm]{Lemma}
\newtheorem{prop}{Proposition}[section]

\theoremstyle{definition}
\newtheorem{defn}{Definition}[section]
\newtheorem{asmp}{Assumption}[section]
\newtheorem{notn}{Notation}[section]

\theoremstyle{remark}
\newtheorem{rmk}{Remark}

\newcommand{\range}{\mathcal{R}}
\newcommand{\acstz}{\mathcal{L}}

\newcommand{\lehm}{\mathcal{M}}

\begin{document}

\newcommand{\fsmpsp}{\left(\Omega, \mathcal{F}, \{\mathcal{F}_t\}, P \right)}
\newcommand{\e}{\mathrm{E}}
\newcommand{\mcal}[1]{\mathcal{#1}}
\newcommand{\dx}{\mathrm{d}}
\newcommand{\ltwo}{{\mathbf{L}^2}}
\newcommand{\linf}{\mathbf{L}^{\infty}}
\newcommand{\lp}{\mathbf{L}^{p}}
\newcommand{\lz}{\mathbf{L}^{0}}
\newcommand{\bb}[1]{\mathbb{#1}}
\newcommand{\norm}[1]{\left\lVert #1 \right\rVert}
\newcommand{\ball}{\left ( \ltwo(S) \right )_{M}} 
\newcommand{\hs}{\mathcal{H}}
\newcommand{\cx}{\mathrm{cx}}
\newcommand{\cc}{\overline{\mathcal{C}}}
\newcommand{\gke}{\epsilon}
\newcommand{\cxl}{\widetilde{\Lambda}}
\newcommand{\cxf}{\tilde{f}}
\newcommand{\hp}{{G}^{\perp}}
\newcommand{\cinth}{\mcal{Z}}
\newcommand{\dqdp}{\frac{\dx Q}{\dx P}}
\newcommand{\tw}{\widetilde{W}}
\newcommand{\ow}{\overline{W}}
\newcommand{\hcl}{\overline{\mcal{G}}}
\newcommand{\acl}{\overline{\mcal{A}}} 
\newcommand{\an}{\mcal{A}_n}
\newcommand{\ltk}{B_k}
\newcommand{\Tau}{\tau}
\newcommand{\zq}{\mcal{Z}_Q}
\newcommand{\rhom}{\rho_{{G}}}
\newcommand{\mk}{\mathcal{K}}
\newcommand{\xiphi}{\mathcal{\pi}}
\newcommand{\leb}{\mathbb{L}}
\newcommand{\bigprob}{\left( P \otimes U_T \otimes \eta \right)}
\newcommand{\medprob}{\left( P \otimes U_T\right)}
\newcommand{\estim}{\mathcal{E}}

\newcommand{\omegafp}{\left(\Omega,\mathcal{F},P\right)}
\newcommand{\calf}{\mathcal{F}}
\newcommand{\omegabmu}{(\widetilde{\Omega} ,\mathcal{B},\mu)}
\newcommand{\otilde}{\Omega\times[T]\times\mathbb{R}}
\newcommand{\mb}{\widehat{\mathcal{F}}\otimes \mathcal{B}(\mathbb{R})}
\newcommand{\omegab}{(\widetilde{\Omega},\mathcal{B})}
\newcommand{\ohat}{\Omega \times [T]}
\newcommand{\fhat}{\widehat{\mathcal{F}}}
\newcommand{\ofnu}{\left(\Omega \times [T], \widehat{\mathcal{F}}, P\otimes U_T\right)}
\newcommand{\varobmu}{(\Theta, \Im, \varrho)}

\newcommand{\lone}{\mathbf{L}^{1}}
\newcommand{\level}{\mathcal{A}}
\newcommand{\power}{\mathcal{T}}
\newcommand{\real}{\mathbb{R}}
\newcommand{\defin}{\stackrel{\triangle}{=}}
\newcommand{\iprod}[1]{\left\langle\; #1\;\right\rangle}
\newcommand{\simp}{\mathcal{S}}
\newcommand{\pset}{\mathcal{P}^{[T]}}
\newcommand{\ftilde}{\widetilde{\mathcal{F}}}
\newcommand{\phitd}{\widetilde{\Phi}}
\newcommand{\dirsum}{\oplus\lone(\calf_t)}
\newcommand{\dirsuminf}{\oplus\linf(\calf_i)}
\newcommand{\rr}{\mathbb{R}}
\newcommand{\br}{\mathbf{B}(\rr)}
\newcommand{\xihat}{\widehat{\xi}}
\newcommand{\h}{Z}
\newcommand{\phiseq}{\mathcal{R}}
\newcommand{\p}{\mathrm{P}}
\newcommand{\C}{\mathcal{C}}
\newcommand{\V}{\mathcal{V}}
\renewcommand{\mcal}[1]{\mathcal{#1}}
\newcommand{\ind}{\mathbb{I}}
\newcommand{\est}{\overline{W.f}}

\newcommand{\var}{V@R}

\theoremstyle{plain}
\newtheorem*{Sauers}{Sauer's Lemma}

\title[Computing Strategies]{Computing Strategies for Achieving Acceptability}
\author{Soumik Pal}
\address{{Department of Statistics, Columbia University}\\ {{1255 Amsterdam Avenue, New York, NY 10027}}}
\curraddr{{{Departments of Mathematics and ORIE}}\\ {{506 Malott Hall Cornell University }}\\ {{Ithaca NY 14853}}}
\email{soumik@math.cornell.edu}
\date{\today}
\thanks{Partially supported by NSF grant 06-01774}

\maketitle

\begin{abstract}
We consider a trader who wants to direct his portfolio towards a set of acceptable wealths given by a convex risk measure. We propose a black-box algorithm, whose inputs are the joint law of stock prices and the convex risk measure, and whose outputs are the numerical values of initial capital requirement and the functional form of a trading strategy to achieve acceptability. We also prove optimality of the obtained capital.   
\end{abstract}


\section{Introduction}

\subsection{Objective}
In this paper, we consider a $T$ period market model, with a single stock and a money market. To model uncertainty in the stock price movements, we consider a probability space $\omegafp$ and a filtration $\calf_0 \subseteq \calf_1 \subseteq \ldots \subseteq \calf_T \subseteq \calf$.
At every time point $t=0,1,2,\ldots,T$, the discounted price of the stock, $S_t$, is assumed to be an integrable random variable measurable with respect to $\calf_t$. 

Next we consider a convex measure of risk. In the following subsection we briefly discuss the definition and significance of such a measure. Here it suffices to define it in the following way. Let $\{Q_i\},\;i=1,\ldots,m$, be a collection of probability measures on the sample space $(\Omega,\calf)$ which are absolutely continuous with respect to $P$, with Radon-Nikod\'{y}m derivatives 
\begin{equation}\label{e:whatisfi}
\{f_i\defin \dx Q_i/\dx P\}.
\end{equation}
We are also given a collection $\{\alpha_i\}$ of real numbers. For every random variable $X \in \cap_{i}\lone(Q_i)$, define
\begin{equation}\label{rho1}
\rho(X)\stackrel{\triangle}{=}\sup_{1\leq i \leq m}\left[\e^{Q_i}(-X) + \alpha_i\right]=\sup_{1\leq i \leq m}\left[-\e(Xf_i) + \alpha_i\right].
\end{equation}
We call such a $\rho$ to be a convex measure of risk.

Let us now introduce an agent who follows a self-financing portfolio by holding $\xi_t$ number of shares in between time periods $t$ and $(t+1)$. Due to the non-anticipative nature of trading, each $\xi_t$ is an $\calf_t$-measurable random variable.  For any choice of initial capital $w_0$, and strategy $(\xi_0,\xi_1,\ldots,\xi_{T-1})$, let $V(w_0,\xi)$ denote the discounted terminal value of the portfolio, i.e., 
\begin{eqnarray}
V(w_0,\xi) &\stackrel{\triangle}{=}& w_0 + W(\xi),\quad \text{where} \label{whatisv}\\
W(\xi) &=& \sum_{t=0}^{T-1} \xi_t(S_{t+1} - S_t). \label{whatisw}
\end{eqnarray}
In this paper we investigate an algorithm to compute a near-minimal $w_0$ and strategy $\xi$, such that $\rho(V(w_0,\xi)) \le 0$. We shall then say that $V(w_0,\xi)$ is \emph{acceptable}.

Our objective is indeed numerical computation, and not just theoretical expressions. We do not impose any restrictions on the law of the price process $S$. However, we do assume the existence of ($\calf_t, P$)-integrable random variables $a_t$ and $b_t$ such that the agent is forced to obey
\begin{equation}\label{restrict}
a_t \leq \xi_t \leq b_t,\;\; \forall\; t=0,1,\ldots,T-1.
\end{equation}
This is often a natural assumptions dictated by trading constraints. In any case, this is crucial for our analysis.

The literature on convex measures of risk is almost silent about computing strategies to achieve acceptability. The primary difficulty being that the terminal conditions on the portfolio are not given by  almost-sure equalities/inequalities. This prevents the use of classical change-of-measure techniques. In this paper, we take an novel computational approach, combining the theory of Uniform Law of Large Numbers with standard Monte-Carlo simulations. 

\subsection{A brief history of the literature}

In recent times, the theory of measures of risk has generated a lot of interest in the mathematical finance literature, partly because it makes a rigorous assessment of risks associated with random financial net worths, and partly because it generalizes No-Arbitrage asset pricing and superhedging ideas in incomplete markets. 

One of the first articles to define and study such measures is the seminal paper \cite{coherent}, which provides a definition and justifies a unified framework for analysis, construction and implementation of \emph{measures of risk}. As the authors point out, these measures of risks, named \emph{coherent} measures, can be used as extra capital requirements, to regulate the risk assumed by market participants, traders, insurance underwriters, as well as to allocate existing capital. The idea is twofold: first to stipulate axioms which define \emph{acceptable} future random net worths, and secondly, to define the measure of risk of an unacceptable position as the minimum extra capital which, invested in a {\lq}pre-specified reference investment instrument{\rq}, makes the future discounted value of the position acceptable. The axioms defining acceptability do not specify a unique measure of risk, instead, they characterize a large class of risk measures. The choice of precisely which measure to use from this class has to be determined from additional economic considerations. 

A significant extension was made by introducing convex measures of risk in \cite{cxrisk}. A similar set-up, as in \cite{coherent}, is considered. However the authors argue that the positive homogeneity of the coherent risk measure is an undue requirement, because the risk of a position might increase in a non-linear way with the size of the position. They suggest to relax the conditions of positive homogeneity and of subadditivity and to require the weaker property of convexity.

In both papers, the basic objects of study are random variables on the set of \emph{states of nature} at a future date, interpreted as possible future (discounted) values of positions or portfolios currently held. A supervisor (e.g. regulator, exchange's clearing firm, or investment manager) decides on a subset of such future outcomes as \emph{acceptable risks}. In other words, they choose a subset $\mcal{A}$ of a suitable set of real functions, $\lz$, on a set $\Omega$, and call it the acceptance set. A measure of risk associated with $\mcal{A}$ is a function $\rho_{\mcal{A}}:\lz \rightarrow \rr$, defined by 
\begin{equation*}
\rho_{\mcal{A}}(X) \defin \inf\{ m \;\vert\; m + X \in \mcal{A}\}.
\end{equation*}
Conversely, for any function $\rho : \lz \rightarrow \rr$, one can define a corresponding acceptance set by $\mcal{A}_{\rho}\defin \left\{ X \in \lz \;\vert \; \rho(X) \le 0\right\}$. Such a function, $\rho$, will be called a \emph{convex measure of risk}, if it satisfies the following axioms:
\smallskip

\noindent$\bullet$ \emph{Translation invariance:} for all $X\in\lz$, and $a\in\rr$, we have $\rho(X + a) = \rho(X) - a$.

\noindent$\bullet$ \emph{Monotonicity:} for all $X$ and $Y$ in $\lz$ with $X \ge Y$, we have $\rho(X) \le \rho(Y)$.

\noindent$\bullet$ \emph{Convexity:} for all $X$ and $Y$ in $\lz$, and all $\lambda\in [0,1]$, we have
\begin{equation}\label{e:cxrisk}
\rho(\lambda X + (1-\lambda) Y) \le \lambda \rho(X) + (1-\lambda)\rho(Y).
\end{equation}
Why these axioms are natural requirements for a measure of risk has been argued in \cite[~Section 2.2]{coherent} and \cite{cxrisk}, and we skip such details. 

The authors of \cite{cxrisk} then prove a representation theorem, similar in spirit to one in \cite{coherent}, which shows that any convex measure of risk on a finite $\Omega$ is of the form
\begin{equation}\label{e:convexrep}
\rho(X) = \sup_{P \in \mcal{P}}\left(\e^{P}[-X] + \alpha(P)\right).
\end{equation}
Here, the set $\mcal{P}$ is the set of all probability measures on $\Omega$. The function $\alpha(\cdot)$ is a certain \emph{penalty function} on $\mcal{P}$ which takes values in $\rr\cup\{-\infty\}$. (Here we stray from the usual convention where $\alpha$ in \eqref{e:convexrep} is replaced by $-\alpha$.) Representation \eqref{e:convexrep} was independently proved by David Heath in \cite{heath}. As before, a convex measure of risk defines an associated acceptance set given by
\begin{equation}\label{e:accconvex}
\mcal{A}_{\rho} = \left\{ X \in \lz\;\vert\; \rho(X) \le 0\right\} = \left\{ X \in \lz\;\vert\; \e^P[X] \ge \alpha(P)\right\}.
\end{equation}
Broad extensions of \eqref{e:convexrep} can be found in \cite{follmer}, all of which exhibit the same structure.

Similar notions as above started appearing simultaneously from very different contexts. In a now well-known paper, \cite{floor}, the authors use the notion of acceptability to present a new approach for positioning, pricing, and hedging in incomplete markets that bridges standard arbitrage pricing and expected utility maximization. Also the theory of \emph{no-good-deal pricing} (NGD), as a pricing technique based on the absence of attractive investment opportunities in equilibrium, was introduced in \cite{cherny}. The term no-good-deal is borrowed from an earlier paper with similar objectives, \cite{cochrane}, where good-deals were defined by high sharp-ratio of returns. The first paper which fully establishes the link between coherent risk measures and the NGD pricing theory is \cite{ngd}, who shows that convex risk measures are essentially equivalent to good-deal bounds. Relations between measures of risk and NGD are further extended by Staum in \cite{staum}, where he proves fundamental theorem of asset pricing for good deal bounds in incomplete markets.

All these diverse motivations can be assimilated by considering what authors of \cite{cxrisk} call \emph{measure of risk in a financial market}. Several authors have recently contributed to the development of this theory, e.g., \cite{barrieu} and \cite{barrieu2}, who establish these risk measures as special cases of \emph{inf-convolution} of risk measures. Consider the setting in the last subsection, in particular, the notations in \eqref{whatisv} and \eqref{whatisw}. The minimum $w_0$ for which $\inf_{\xi}\rho(V(w_0,\xi))$ is non-positive can be thought of as a price one has to pay today for achieving acceptability in future. As is shown in \cite{cxrisk}, for any random variable $Z$, one can choose the penalty function suitably such that the minimum $w_0$ is the market measure of risk of $Z$. This duality between price and risk measures is also seen in NGD pricing. If a strategy $\xi$ exists which achieves the infimum above, then, it can be thought of as a hedging strategy in the NGD setting. In any case, it can be thought as a strategy to achieve acceptability in the future, starting from a currently non-acceptable portfolio.

Our risk measure, $\rho$ in \eqref{rho1}, is clearly convex. We restrict ourselves to finite sets $\{Q_1,\ldots,Q_m\}$ to aid computation. This can be interesting either in its own right, or as an approximation to the general case. The assumption $Q_i \ll P$ is implied by the natural requirement: $\rho(X)=\rho(Y)$ if $P(X=Y)=1$ (see, e.g., \cite{follmer}).

\subsection{Summary and organization} 

We propose our main result in the following section. First, we suppose that for a given $w_0$, the set of strategies $\xi$ which satisfy \eqref{restrict}, and for which $\rho(V(w_0,\xi)) \le 0$ is non-empty. Then, Proposition \ref{main} proves that the intersection of this set with a specific, much smaller family of strategies is also non-empty. This smaller set of strategies is indexed by a finite-dimensional space, and has nice combinatorial properties. This allows us to use the theory of Uniform Law of Large Numbers (ULLN), and devise a Monte-Carlo scheme to numerically compute a near-minimum $w_0$ and a corresponding strategy $\xi$ to have $\rho(V(w_0,\xi))$ non-positive. In Section \ref{compu}, we describe the method, and give precise error bounds on such approximations. In Section \ref{examples}, we consider a natural example in which stock price follows discrete geometric Brownian motion, and show how our method leads to numerical values of both near-optimal capital and strategy to achieve acceptability.

\subsection{Acknowledgments}
I thank Prof. Peter Bank for suggesting the particular example in Section \ref{examples}.

\section{Main results} Recall that $m$ refers to the number of probability measures in the representation of $\rho$ in \eqref{rho1}. Let $\acstz$ be the collection of adapted processes $\xi=(\xi_0,\ldots,\xi_{T-1})$, which satisfy \eqref{restrict}. Define the following set:
\begin{equation}\label{whatisr}
\range\defin \left\{\left(\e^{Q_1}(W(\xi)),\ldots,\e^{Q_m}(W(\xi))\right),\quad \xi \in \acstz \right\}.
\end{equation}
For any $k\in\mathbb{N}$ and any $x \in \real^{k}$, define the \emph{upper quantant} of $x$, denoted by $Q_x$, as the set
\begin{equation}\label{lquant}
Q_x=\{y\in\real^{k}: y_j\geq x_j \;\mbox{for}\; j=1,2,\ldots,k\}. 
\end{equation}
The dimension is suppressed in the notation for $Q_x$, since it is obvious from the dimension of $x$.

\begin{prop}\label{main}

Fix a $w_0\in\rr$. Let $z_0=(\alpha_1-w_0,\ldots,\alpha_m-w_0)$. Assume that the convex set $Q_{z_0}\cap\range$ has a non-empty relative interior. 

For every $1\le i\le m$, define the adapted sequence of random variables
\begin{equation}\label{whatisvt}
v_t(f_i) \defin (b_t-a_t)\e\left[(S_{t+1}-S_{t})f_i\;\Big\vert\;\calf_t\right],\;\; t=0,1,\ldots,T-1.
\end{equation}
For every $\mathbf{r}\in \real^{m}$, consider the following weighted sum process
\begin{equation}\label{whatislambdat}
\lambda_t(\mathbf{r}) \defin \sum_{i=1}^m r_iv_t(f_i), \;\;t=0,\ldots,T-1.  
\end{equation}
Now, let $\eta$ be any continuous probability distribution function on the real line with finite first moment. Then, there exists a vector $\mathbf{r^*}\in\real^{m}$, such that the $\{\calf_t\}$-adapted process
\begin{equation}\label{eq:compu}
\xi^*_t(\omega) \defin (b_t - a_t)\eta(\lambda_t(\mathbf{r^*})) + a_t, \;\;t=0,\ldots,T-1,
\end{equation}
satisfies~(\ref{restrict}) and $\rho(W(\xi^*)) \le w_0$.
\end{prop}

\begin{rmk}
Note that assuming $Q_{z_0}\cap\range$ being non-empty is equivalent to assuming the existence of a strategy $\xi$ such that $\rho(w_0 +W(\xi)) \le 0$. In the above proposition we assume a bit more than that.

\end{rmk}

The proof of this result will follow after we have introduced some notations. Let $[T]$ denote the set $\{0,1,\ldots,T-1\}$. Enlarge the original sample space by considering
\begin{equation}\label{otilde}
\Omega \times [T] = \Omega \times \{0,1,\ldots,T-1\}.
\end{equation}
Let $\pset$ be the power set of the finite collection $\{0,1,\ldots,T-1\}$ and let $\calf \otimes \pset$ denote the product $\sigma$-algebra between $\calf$ and $\pset$. Extract a sub $\sigma$-algebra $\fhat$ by defining
\begin{equation}\label{whatismb}
\fhat \defin \left\{ A\in\calf\otimes\pset\;\Big\vert\; \{ \omega\; :\; (\omega, t)\in A \} \in \calf_t, \;\; \forall \; t=0,1,\ldots, T-1\right\}.
\end{equation}
That $\fhat$ is a valid $\sigma$-algebra is straightforward to verify. Finally, let $U_T$ denote the discrete uniform measure on $[T]$, and consider the product measure $P\otimes U_T$ on the $\sigma$-algebra $\fhat$. This gives us a probability space $\ofnu$. The advantages of considering the above probability space is the following trivial lemma.

\begin{lem}\label{bijections}
A process $(h_0,h_1,\ldots,h_{T-1})$ is adapted with respect to $(\Omega, \calf)$ if and only if the random variable $H(\omega,t)=h_t(\omega)$ is measurable with respect to the enlarged space $(\Omega\times [T], \fhat)$. 
\end{lem}

\begin{proof} Follows from the definition of $\fhat$.
\end{proof}

For all sequence $\{\xi_t\}$ that satisfy~(\ref{restrict}) (i.e. $\xi\in\acstz$), let us make a change of variable $\phi=\xiphi(\xi)$, where
\begin{equation}\label{phixi}
\phi_t=\xiphi(\xi)_t\defin(\xi_t - a_t)/(b_t - a_t),\;\; t=0,\ldots,T-1,
\end{equation}
then, each $\phi_t$ is $\calf_t$-measurable and $P(0\le \phi_t \le 1)=1$.
 
Now, the discounted terminal value of the portfolio in~(\ref{whatisw}) can be expressed in terms of the $\phi = \xiphi(\xi)$ as
\begin{equation}
\begin{split}
W(\xi) = W\circ \xiphi^{-1}(\phi)&= \sum_{t=0}^{T-1} \left(S_{t+1} - S_t\right)\left[ (b_t-a_t)\phi_t + a_t\right]\\
 &= \sum_{t=0}^{T-1} \left[(b_t - a_t)\left(S_{t+1} - S_t\right)\phi_t + a_t (S_{t+1} - S_t)\right] 
\end{split}
\end{equation}
Thus, for any suitably integrable $f$ defined on $\omegafp$, one can write
\begin{eqnarray}\label{winphi}
\int W(\xi) f\dx P = \e\left(W(\xi)f\right) &=& \sum_{t=0}^{T-1} \e\Big(\left[(b_t - a_t)\left(S_{t+1} - S_t\right)\phi_t + a_t (S_{t+1} - S_t)\right] f\Big)\nonumber\\
 &=&  \sum_{t=0}^{n-1} \e\left[ (b_t - a_t)\left(S_{t+1} - S_t\right)\phi_t f\right] +  \sum_{t=0}^{n-1} \e \left[ a_t (S_{t+1} - S_t) f \right] \nonumber\\
 &=& \sum_{t=0}^{T-1} \e\left[ v_t(f) \phi_t\right] + c(f),
\end{eqnarray}
where, we have named
\begin{subequations}
\begin{equation}\label{whatisvtf}
v_t(f)(\omega) \defin (b_t - a_t)\e\left[ \left(S_{t+1} - S_t\right)f \Big\vert \calf_t\right](\omega), \quad \text{and}
\end{equation}
\begin{equation}\label{whatiscf}
c(f) = \e\left[f\sum_{t=0}^{T-1} a_t (S_{t+1} - S_t) \right].
\end{equation}
For $ t \in [T]$, if we now look at $\phi$ and $v$ as functions of two arguments $(\omega, t)$, i.e.,
\begin{equation}\label{enlarge}
\phi(\omega,t) \defin \phi_t(\omega),\;\; v(f)(\omega,t) \defin v_t(f)(\omega),\;\; \omega\in\Omega,
\end{equation}
\end{subequations}
then, by Lemma~\ref{bijections}, both $\phi$ and $v$ are $\fhat$-measurable functions on $\ohat$. Moreover, $v(f)$ is $P\otimes U_T$-integrable and $P\otimes U_T(\{0\leq \phi \leq 1\})=1$. Thus, from~(\ref{winphi}), we can write
\begin{equation}\label{linform}
\begin{split}
\int W\circ\pi^{-1}(\phi) f \dx P - c(f) = \sum_{t=0}^{T-1} \int v_t(f)\phi_t\dx P = T\int_{\Omega\times [T]} \phi v(f) \dx \left(P\otimes U_T\right).
\end{split}
\end{equation}

\begin{proof}[Proof of Proposition~\ref{main}] Let $\eta$ be any continuous probability distribution function on the real line with finite first moment. Consider the probability space $(\rr, \mcal{B}(\rr), \eta)$, where $\mcal{B}(\rr)$ is the Borel $\sigma$-algebra on $\rr$. Consider the following product space
\begin{equation}\label{obmu}
\otilde,\quad \mb ,\quad P \otimes U_T\otimes \eta.
\end{equation}
Let us recall here that $\ohat$, and $\fhat$ are defined in~(\ref{otilde}),~(\ref{whatismb}), and $U_T$ is the discrete uniform measure on the set $[T]=\{0,1,2,\ldots, T-1\}$. Let $\h$ be a measurable map from this product space to $\rr$, given by
\[
\h(\omega, t, x) = x, \qquad \omega \in \Omega,\; t \in [T],\; x \in \rr.
\]
Clearly, $\h$ has distribution $\eta$, independent of the $\sigma$-algebra $\fhat$.
 
Consider the functions $f_i=\dx Q_i/\dx P$ appearing in \eqref{rho1}, and define the following functions in $\lone(P\otimes U_T \otimes \eta)$:
\begin{equation}\label{whatisg}
 g_i(\omega,t,x) \defin v(f_i)(\omega,t),\;\; 1\le i\le m,\quad g_{m+1} \defin -\h,
\end{equation}
where the function $v$ is defined in \eqref{whatisvtf} and ~(\ref{enlarge}). Also define the constants 
\[
\gamma_i \defin \left(\alpha_i - w_0 - c(f_i)\right)/T, \;\; i=1,2\ldots,m.
\]
The function $c$ is defined above in \eqref{whatiscf}.

\noindent $\bullet$ Define $\Phi$ to be the convex collection of all $\mb$-measurable functions $\phi$ such that $P\otimes U_T \otimes \eta(0\le \phi\le 1) = 1$. Let $\lehm$ denote the set of points
\[
\lehm\defin \left\{\left(\int\phi g_1,\ldots,\int\phi g_{m+1}\right),\quad \phi \in \Phi\right\},
\]
where the integrations are with respect to $P\otimes U_T \otimes \eta$. Recall the assumption in the statement of the proposition that $Q_{z_0}\cap\range$ has a non-empty relative interior. Since every strategy $\xi \in \acstz$ corresponds to a $\phi$ by the linear mapping $\pi$ defined in \eqref{phixi}, it follows that there is a point $(q_1,\ldots,q_m)$ which is an interior point of $\lehm \cap Q_{\gamma}$.

We look at the following maximization problem: find the maximizer of 
\[
\int_{\Omega\times[T]\times\rr}\phi g_{m+1}\dx\left(P\otimes U_T\otimes \eta\right)
\]
among all $\phi \in \mcal{A} \subseteq \Phi$, where $\mcal{A}$ is defined by
\begin{equation}\label{e:ainrevision}
\mcal{A} \defin \left\{ \phi \in \Phi \;\bigg\vert\; \int_{\Omega\times[T]\times\rr}\phi g_{i}\dx\left(P\otimes U_T\otimes \eta\right) = q_i \right\}.
\end{equation}

We use Theorem 5 on page 96 of \cite{lehmann}. Part (iv) of this Theorem guarantees the existence of a solution $\phi^*$ of the above maximization problem which is of the form
\begin{equation}\label{phi*}
\phi^* = \left\{ \begin{array}{ll}
1, & \text{if} \;\;\sum_{i=1}^{m} r_iv(f_i) + g_{m+1} > 0,\\
\\
0, & \text{if} \;\;\sum_{i=1}^{m} r_iv(f_i) + g_{m+1} < 0,
\end{array}\right.
\end{equation}
for some $(r_1,\ldots,r_{m})\in \rr^{m}$.
Recall that $g_{m+1}=-Z$, and that from the definition of the function $v$ in~(\ref{enlarge}), it is clear that each $v(f_i)$ is independent of $Z$. Thus
\begin{equation*}
\begin{split}
\bigprob\left( \sum_{i=1}^{m} r_iv(f_i) = Z\right) =
\int \eta\left(\h =\sum_{i=1}^{m} r_iv(f_i)\;\Big\vert\;\;\fhat\right)\dx \medprob =0,
\end{split}
\end{equation*}
the integrand is zero being the consequence of the continuity of $\eta$. Thus, the solution in~(\ref{phi*}) is actually 
\begin{equation}
\phi^* = \begin{cases}
1 & \text{if} \;\;\sum_{i=1}^{m} r_iv(f_i) > Z \\
0 & \text{otherwise}.
\end{cases}
\end{equation}
Now from the constraint $\phi^* \in \mcal{A}$, we also get that $\int \phi^* g_i\dx\bigprob = q_i \ge \gamma_i,\;\; i=1,2,\ldots,m$. In other words, $\int \phi^* v(f_i)\dx \bigprob \ge \gamma_i$, or, by expanding $\gamma_i$, equivalently
\begin{equation}\label{link}
T\int_{\ohat}\e\left( \phi^* \;\Big\vert\;\fhat\right) v(f_i)\dx \left(P\otimes U_T\right) + c(f_i) \geq \alpha_i - w_0,\;\; i=1,\ldots,m.
\end{equation}
We have used Fubini above, where
\begin{eqnarray*}\label{whatispsi}
\e\left( \phi^* \;\Big\vert\;\fhat\right)(\omega,t) &=& \bigprob\left(  \sum_{i=1}^{m} r_iv(f_i) - \h > 0\;\Big\vert\;\fhat\right)(\omega, t)\\
&=& \eta\left(\sum_{i=1}^{m}r_iv_t(f_i)(\omega)\right) = \eta\left(\lambda_t(\omega)\right),
\end{eqnarray*}
where the $\{\calf_t\}$-adapted process $\{\lambda_t\}$ is defined as in~(\ref{whatislambdat}). Thus, if we let
\[
\xi^*_t\defin (b_t-a_t)\eta(\lambda_t) + a_t,\;\;t=0,\ldots,T-1,
\]
then, $\xi^*=\pi^{-1}(\e[ \phi^* \vert\fhat])$ in the notation of \eqref{phixi}. Thus, by~(\ref{link}) and~(\ref{linform}), we conclude that 
\[
\int W(\xi^*)f_i\dx P \ge \alpha_i - w_0,\;\; i=1,\ldots,m,
\]
or in other words, $\rho(w_0+W(\xi^*)) \le 0$. This proves the proposition.\end{proof}

\section{Computations}\label{compu}

For every ${s}=(s_1,\ldots,s_{m}) \in \real^{m}$, recall from Proposition~\ref{main}, the $\calf_t$-adapted process 
\begin{equation}\label{e:whatislambdats}
\lambda_t(s) \defin \sum_{i=1}^{m}s_iv_t(f_i),\quad t=0,1,\ldots, T-1,
\end{equation}
and the derived process
\begin{equation}\label{remindxit}
\xi_t(s) \defin (b_t-a_t)\eta(\lambda_t(s)) + a_t,\quad t=0,\ldots,T-1. 
\end{equation}
For suitable $w_0$, Proposition~\ref{main} proves the existence of an $s \in \real^{m}$ via which the process $\xi(s)$ satisfies $\rho(w_0 + W(\xi(s))) \le 0$, or equivalently, by translation invariance, $\rho(W(\xi(s))) \le w_0$. 

What we shall do now is like a partial converse. Suppose we can compute $\rho(W(\xi(s)))$ for every $s\in\real^m$. Then we can define $w_0$ by
\[
w_0 := \inf_{s} \rho(W(\xi(s))).
\]
If $s^*$ achieves the above infimum, then, clearly $\rho(w_0+ W(\xi(s^*)) \le 0$, and $w_0$ is near minimal by Proposition~\ref{main}.

The above procedure would work if we could theoretically compute $\rho(W(\xi(s)))$ for every $s \in \real^{m}$. This is often impossible. 
However, for any fixed $s$, we can estimate $\rho(W(\xi(s)))$ by Monte-Carlo simulations upto any desired level of accuracy. We show in this section that it is possible to do a Monte-Carlo simulation to simultaneously approximate $\rho(W(\xi(s)))$ for \emph{every} $s\in\real^{m}$ with a uniform error bound. The feasibility of our claim depends on the theory of {\textit{Uniform Law of Large Tumbers}} and the related concept of Vapnik-\u{C}ervonenkis dimension which is a combinatorial property of the particular structure of $\{\xi_t\}$ in \eqref{remindxit}. This theory is well-developed and we cherry-pick only the necessary results for our purpose. These have been stated in the appendix. Further references have also been provided for the interested reader.  

Once we have our estimation procedure ready, we can construct a finite mesh $\mathbb{G}$ within $\rr^{m}$ and approximate the value of $\rho(W(\xi(r)))$, by (say) $\hat{\rho}(r)$, for every $r \in \mathbb{G}$. Note that the error in approximation in our Monte-Carlo procedure does not depend on the size of the grid, and we can make it as large and fine as we want. For that fine mesh $\mathbb{G}$, let $r^*$ be a grid point which attains $\hat{\rho}(r^*) = \min_{r \in \mathbb{G}} \hat{\rho}(r)$.

Let $w_0^* = \hat{\rho}(r^*)$. Then, as we describe below, given any $\epsilon, \delta > 0$, with a very high probability greater than $(1-\delta)$, the choice of $(w_0^*, \xi(r^*))$ satisfies
\[
\rho(w_0^* + W(\xi(r^*))) \le \epsilon.
\]
This gives a near-minimal initial capital for the problem of finding $(w_0, \xi)$ which satisfies \eqref{restrict} and $\rho(w_0 + W(\xi)) \le \epsilon$. 
\bigskip

Central to computing $\rho(W(\xi(s)))$, for any $s \in \simp_{m+1}$, is to compute $\e(W(\xi(s))\cdot f_i)$ for every $f_i$ that defines $\rho$. Now, from equation \eqref{winphi}, we can write
\begin{eqnarray}
\e(W(\xi(s)) f_i) &=& \sum_{t=0}^{T-1} \int_{\Omega} v_t(f_i)\eta(\lambda_t(s)) \dx P + c(f_i)\nonumber\\
&=& T \int_{\Omega\times [T]}\eta(\lambda(s)) v(f_i) \dx \medprob + c(f_i).\nonumber\\
&=& T \int_{\otilde}\ind\left\{ \lambda(s) - Z > 0\right\} v(f_i) \dx \bigprob + c(f_i).\label{asanind}
\end{eqnarray}
Here, as in the last section, $Z$ is a random variable with law $\eta$ independent of $\fhat$, and $\mathbb{I}\{\cdot\}$ denotes the indicator of an event.

We would now like to do a change of measure in \eqref{asanind} above with $v(f_i)$ as the {\lq}Radon-Nikod\'{y}m{\rq} derivative. This is not possibly directly, since $v(f_i)$ is not necessarily positive. However, we can work separately with $v^+(f_i)=\max(v(f_i), 0)$ and $v^-(f_i)=\max(-v(f_i),0)$, which denote the positive and the negative parts respectively. Hence, one obtains
\begin{eqnarray}\label{compu:ind}
\e(W(\xi(s)) f_i) - c(f_i)
&=& T \int_{\otilde}\ind\{\lambda(s) - Z > 0\}\; v^{+}(f_i)\; \dx \bigprob \nonumber\\
& &- T \int_{\otilde}\ind\{\lambda(s) - Z > 0\}\; v^{-}(f_i)\; \dx \bigprob\nonumber\\
&=& d_i^{+}\cdot\left(\mu_i^{+}\otimes\eta\right)\{\lambda(s) - Z > 0\} \nonumber\\
& & - d_i^{-}\cdot\left(\mu_i^{-}\otimes\eta\right)\{\lambda(s) - Z > 0\}.
\end{eqnarray}
Here we have introduced several probability measures on $(\Omega\times[T], \fhat)$, defined by their corresponding \emph{unnormalized} Radon-Nikod\'{y}m derivatives:
\begin{subequations}
\begin{equation}\label{Qpm}
{\dx \mu_i^{+}}/{\dx\medprob} \propto v^{+}(f_i), \qquad {\dx \mu_i^{-}}/{\dx \medprob} \propto v^{-}(f_i),
\end{equation}
and the corresponding normalizing constants (multiplied by $T$):
\begin{equation}\label{alphapm}
d_i^{+} \defin \sum_{t=0}^{T-1} \e[v^{+}_t(f_i)], \quad d_i^{-} \defin \sum_{t=0}^{T-1} \e[v^{-}_t(f_i)],\quad i=1,2,\ldots,m.
\end{equation}
\end{subequations}
If any of the constants in \eqref{alphapm} is zero, the corresponding measure becomes the zero measure and can be dropped from our analysis. For efficiency in computation we would like to keep track of the number of non-zero measures above by defining
\begin{equation}\label{howmanyzero}
\aleph \defin \sum_{i=1}^m \left( 1_{\{ d_i^+ > 0\}} + 1_{\{ d_i^- > 0\}}\right). 
\end{equation}
\bigskip

\begin{asmp}\label{compu:main}
Throughout the rest of this section, we shall assume that 
\begin{enumerate}
\item one can generate samples from the joint distribution of $(S_0,S_1, \ldots, S_T)$, 
\item the random variables $v_t(f_i)$ (thus also $\lambda_t$) can be evaluated given the values of $(S_0,\ldots, S_T)$, and 
\item the constants $c(f_i)$, $d^+_i$ and $d^-_i$ can be evaluated for every $1\le i \le m$.
\end{enumerate}
\end{asmp}

Now, by \eqref{compu:ind}, evaluating $\e(W(\xi(s))f_i)$ boils down to evaluating the following two probabilities 
\begin{equation}\label{e:probtoestim}
\left(\mu_i^{+}\otimes\eta\right)\{\lambda(s) - Z > 0\},\quad \text{and}\quad \left(\mu_i^{-}\otimes\eta\right)\{\lambda(s) - Z > 0\},\quad s\in \rr^m.
\end{equation}
Instead, we use the Vapnik-\u{C}ervonenkis theory, described in Subsection \ref{subsec:ulln} in the Appendix, to set up a Monte-Carlo scheme to estimate them for all $s\in\real^m$ with uniform precision. The key to this is to observe the trivial equality 
\begin{equation}\label{e:vccompute}
\{ \lambda(s) - Z > 0 \} = \left\{ \sum_{j=1}^{m} s_j v(f_j) - Z > 0 \right\}
\end{equation}
and apply Dudley's Theorem, Theorem \ref{dud} in the Appendix, with $X = \Omega \times [T] \times \rr$ and the vector space $G$ to be linear space spanned by $Z$ and $v(f_j)$, $j=1,2,\ldots,m$. Thus we  infer that the collection of sets  
\begin{equation}\label{whatisd}
\Biggl\{ \biggl\{ \tilde{\omega}\in\otilde\;:\;\sum_{j=1}^{m} r_j v(f_j)(\tilde{\omega}) + r_{m+1} Z(\tilde{\omega}) > 0 \biggr\},\; r \in \rr^{m+1}\Biggr\},
\end{equation} 
has a VC dimension not more than $(m+1)$. From \eqref{e:vccompute}, the collection of sets
\[
\left\{\{ \lambda(s) - Z > 0 \}, \quad s \in \real^{m}\right\}
\]
is contained in \eqref{whatisd}, and hence also has a VC-dimension not more than $(m+1)$.
It is hence possible to estimate the probabilities in 
\eqref{e:probtoestim}, uniformly for all $s \in \real^{m}$, by drawing independent samples from distributions $\mu_i^{+}\otimes\eta$ and $\mu^{-}_i\otimes\eta$. 
\bigskip

Our aim now would be to apply Theorem \ref{combovc}. We first have to choose two positive parameters, $\epsilon$ and $\delta$, determining the precision of our estimates. Now, for every $i=1,2,\ldots,m$, choose $\kappa_i^+$ such that
\begin{equation}\label{e:nosamples}
4(\kappa^{+}_i)^{2(m+1)}\exp\left( -2\kappa^{+}_i\left(\frac{\epsilon}{d^{+}_i}\right)^2 + 4\left(\frac{\epsilon}{d^{+}_i}\right) + 4\left(\frac{\epsilon}{d^{+}_i}\right)^2\right) \le \delta.
\end{equation}
Generate $\kappa_i^+$ many iid samples $\{(\omega_j, t_j, z_j) \in \Omega\times [T] \times \rr$, $j=1,2,\ldots, \kappa_i^+\}$, from the joint distribution $\mu^+_i\otimes\eta$. 

\begin{rmk}
It is fairly standard to generate samples from measures $\mu^+_i$, defined through their unnormalized densities given in \eqref{Qpm}. We can either directly identify the distribution, as we do in the next section. Or, under the assumption that one can generate perfect samples from the underlying distribution $\medprob$, one can use any of the standard Markov Chain algorithms, from the simple rejection sampling, to the general Metropolis-Hastings algorithm to generate samples from $\mu^{+}_i$. Several books, e.g. \cite[~Chap. 11]{BDA}, describe the details of all these algorithms. 
\end{rmk}

Let $\estim_i^{+}(\cdot)$ denote the empirical estimates of probabilities by the sample frequency. For example, for any $s\in \real_{m}$, we have
\begin{equation}\label{e:estim}
\estim_i^+\{\lambda_s - Z > 0\} = \frac{1}{\kappa_i^+}\sum_{j=1}^{\kappa_i^+} \mathbb{I}{\{\lambda_{t_j}(s)(\omega_j) - z_j > 0\}}.
\end{equation}
We can now apply \eqref{ulln} from Theorem \ref{combovc} to claim that under the joint distribution of all the $\kappa_i^{+}$ many samples drawn
\begin{equation}\label{approx+}
\text{Prob}\left\{ \sup_{s\in\real^{m}} d^{+}_i\Large\lvert \;\estim_i^+\{\lambda_s - Z > 0\}-({\mu_i^+\otimes\eta})\{\lambda_s - Z > 0\}\; \Large\rvert > \epsilon \right\} \le \delta, \;\forall i.
\end{equation}
Exactly in the same way, one can replace the $\mu^+_i$ by $\mu_i^{-}$ above, compute $\kappa_i^-$ by 
\begin{equation}\label{e:nosamples-}
4(\kappa^{-}_i)^{2(m+1)}\exp\left( -2\kappa^{-}_i\left(\frac{\epsilon}{d^{-}_i}\right)^2 + 4\left(\frac{\epsilon}{d^{-}_i}\right) + 4\left(\frac{\epsilon}{d^{-}_i}\right)^2\right) \le \delta,
\end{equation}
and obtain estimates $\estim_i^{-}$, analogous to \eqref{e:estim}, which satisfies
\begin{equation}\label{approx-}
\text{Prob}\left\{ \sup_{s\in\simp_{m+1}} d^{-}_i\Large\lvert \;\estim_i^-\{\lambda_s - Z > 0\}-({\mu_i^-\otimes\eta})\{\lambda_s - Z > 0\}\; \Large\rvert > \epsilon \right\} \le \delta, \;\forall i.
\end{equation}
From \eqref{approx+} and \eqref{approx-}, it follows, by using \eqref{compu:ind}, that one can estimate the quantity $\e(-W(\xi(s))f_i) + \alpha_i$ by 
\begin{equation}\label{whatisdi}
\mcal{D}_i(s) \defin -d_i^{+}\estim_i^{+}\{\lambda(s) - Z > 0\}  + d_i^{-}\estim_i^-\{\lambda(s) - Z > 0\} - c(f_i) + \alpha_i. 
\end{equation}
Since $\rho(W(\xi(s)))= \sup_{1\le i \le m} \{ \e(-W(\xi(s))f_i) + \alpha_i\}$, it follows that a good estimate of $\rho(W(\xi(s)))$ is 
\[
\hat{\rho}(s)\defin\sup_i\mcal{D}_i(s).
\]
We can sum-up this approximation by a simple 
union bound using \eqref{approx+} and \eqref{approx-} as follows.
\medskip

Under the joint distribution of all the $\{\kappa_i^+, \kappa_i^-\}_{1\le i \le m}$ samples drawn from the distributions $\{ \mu_i^+\otimes\eta,\; \mu_i^-\otimes\eta\}_{1\le i \le m}$, one has
\[
\text{Prob}\left\{ \sup_{s\in\real^{m}} \;\lvert \hat{\rho}(s) - \rho(W(\xi(s)))\;\rvert \ge \epsilon \right\} \ge 1 - \aleph \delta.
\]
Here, the number $\aleph$ ($\le 2m$) is described in \eqref{howmanyzero}. We use the number $\aleph$ and not the crude bound $2m$ to bring more efficiency in our estimate.  

Now that we have estimated $\rho(W(\xi(s)))$ for every $s \in \real^{m}$ with uniform precision, we can carry out the grid searching procedure described at the beginning of this section. We minimize $\hat{\rho}(s)$ over the grid nodes (say $\mathbb{G}$) to obtain
\[
w_0^* \defin \inf_{s\in\mathbb{G}}\hat{\rho}(s)= \hat{\rho}(s^*).
\]
Then, with a probability more than $(1 - \aleph \delta)$, we have
\[
\rho(W(\xi(s^*)) \le \hat{\rho}(s^*) + \epsilon \le w_0^* + \epsilon.
\]
In other words, with a high probability of being correct, we get
\[
\rho(w_0^* + W(\xi(s^*)) \le \epsilon.
\]
Thus one obtains a near-optimal pair $(w_0, \xi)$ which satisfies \eqref{restrict} and $\rho(w_0 + W(\xi))$ is \emph{almost} non-positive. The next section displays the entire method through an explicit example.

\section{Examples}\label{examples} The previous theory is now applied to an explicit example where stock prices follow geometric Brownian motion, but observed only at finitely many time points.

We consider $T=3$ and $\Omega= \rr^T$, the $\sigma$-algebra $\calf_t$ being generated by the first $t$ co-ordinates of $\omega \in \Omega$. We take $\calf_0$ to be the trivial $\sigma$-algebra $\{\emptyset, \Omega\}$. Take $P$ to be the product probability measure of $T$ many independent Normal distributions with mean zero and variance one. In other words, we consider random variables $(Z_1,Z_2, \ldots, Z_T)$ such that each $Z_i$ is independent and identically distributed as $N(0,1)$. The discounted stock price movement, under $P$, is described by
\begin{equation}\label{exmp:model}
S_0=4, \quad S_{t+1}={S_t}\exp\left[ -\frac{1}{2} + Z_{t+1} \right], \; t=0,1,\ldots,T-1.
\end{equation}
In other words, we have
\begin{equation}\label{modelrestate}
S_t = S_0 \exp\left[ \sum_{i=1}^tZ_k - \frac{t}{2} \right], \quad t=1,2,\ldots, T.
\end{equation}
However, the investor is not entirely certain of his modeling assumptions, and so considers other scenarios $Q_1$ and $Q_2$, where $Q_1$ and $Q_2$ are two probability measures defined on $(\Omega,\calf_T)$ by
\[
\begin{split}
\text{under}\; Q_1, \qquad& Z_1,\ldots, Z_T \stackrel{iid}{\sim} N(1,1),\\
\text{under}\; Q_2, \qquad& Z_1,\ldots, Z_T \stackrel{iid}{\sim} N(-1,1).
\end{split}
\] 
For convenience we also introduce $Q_3 = P$.

\begin{rmk}
Note, from \eqref{exmp:model}, the effect of changing measure on the stock price movements. For $Q_1$, the geometric Brownian motion gets a positive drift, for $Q_2$ it gets a negative drift, while $Q_3$ is the same as $P$, where stock prices are a martingale.
\end{rmk}

Assume that various constraints dictate that his trading strategy is bounded between zero and one throughout, i.e., in the notation of \eqref{restrict}, we have
\[
a_t \equiv 0, \qquad b_t \equiv 1, \quad \text{for all}\quad 0\le t\le T-1.
\]
Now, the investor sets to do the following: if the conditions are favorable, and the stock prices tend to go up under $Q_1$, he wants a large lower bound $\mathrm{e}^4$ for his expected terminal wealth. On the other hand, if the stock prices tend to go down, under $Q_2$, he sets a lower bound for his expected losses, by setting that his final expected wealth should be more than $\mathrm{e}^{-1}$. He has at least \$$0.2$ to invest, and would like to know an optimal initial capital, and a trading strategy to achieve his goals. 

This requires us to define a measure of risk $\rho$: if $X$ is measurable with respect to $\calf_T$, then 
\[
\rho(X) \defin \max_{i=1,2,3}[ \e^{Q_i}(-X) + \alpha_i ],\qquad m=3,
\]
where
\[
\alpha_1 = \mathrm{e}^4,\quad \alpha_2 = \mathrm{e}^{-1},\quad \alpha_3 = 0.2.
\]
Then, we would like to compute a near-optimal pair $(w_0, \xi)$ of initial capital $w_0$ and $0 \le \xi_t \le 1$, for all $0\le t\le T-1$, such that 
\[
\rho(w_0 + W(\xi)) \le 0 \quad \Leftrightarrow \quad w_0 + \e^{Q_i}[W(\xi)] \; \ge \alpha_i, \quad i=1,2,3.
\]
\medskip

The first step will be to compute the functions $f_1, f_2,$ and $f_3$. They are straightforward since
\begin{equation}\label{computefi}
\begin{split}
f_1(z_1,\ldots,z_k) &= \dx Q_1/ \dx P = \exp\left[ \sum_{k=1}^T z_k - T/2\right]\\
f_2(z_1,\ldots,z_k) &= \dx Q_2/ \dx P = \exp\left[ -\sum_{k=1}^T z_k - T/2\right]\\
f_3(z_1,\ldots,z_k) &= \dx Q_3/ \dx P \equiv 1.
\end{split}
\end{equation}
We can now compute the functions $v_t(f_i)$. These are given by
\begin{equation}
\begin{split}
v_t(f_1) &= \e\left[ f_1 (S_{t+1} - S_t) \;\vert \;\calf_t\right]\\
 &= S_t\e\left[ f_1 \left(\exp(Z_{t+1}-1/2) - 1\right) \;\vert \;\calf_t\right],\quad \text{from}~\eqref{exmp:model},\\
 &= S_t\exp\left(\sum_{k=1}^t Z_k\right)\e\left( \exp\left\{\sum_{k=t+1}^T Z_k - T/2\right\}\left[ \exp(Z_{t+1}-1/2) - 1 \right]\right),
\end{split}
\end{equation} 
where the last equality is due to \eqref{computefi} and the independence of $\{Z_i\}$. Recall that if $Z$ follows $N(0,1)$, then
$\e\left[ \exp(\sigma Z) \right] = \exp( \sigma^2/ 2), \; \sigma \in \rr$. Thus, for $\mathbf{z}=(z_1,z_2,\ldots,z_T) \in \Omega$, a straightforward computation leads to
\begin{equation}\label{exmp:whatisvtf1}
\begin{split}
v_t(f_1)(\mathbf{z}) &= S_t \exp\left[ \sum_{k=1}^t z_k\right]\left\{ \exp\left( 1-\frac{t}{2}\right) - \exp\left(-\frac{t}{2}\right)\right\}\\
&= 4(\mathrm{e}-1)\exp\left\{2\sum_1^t z_k - t\right\}, \quad \text{by}\; \eqref{modelrestate}.
\end{split}
\end{equation}
In particular, we have $\e(v_t(f_1)) = 4(\mathrm{e}-1)\e\left[\exp\left(2\sum_{k=1}^t Z_k -t \right)  \right]= 6.87 \mathrm{e}^{t}$.
 
Similarly, we compute
\begin{eqnarray}\label{exmp:whatisvtf2}
v_t(f_2) &=& \e\left[ f_2 (S_{t+1} - S_t) \;\vert \;\calf_t\right]\nonumber\\
 &=& S_t\e\left[ f_2 \left(\exp(Z_{t+1}-1/2) - 1\right) \;\vert \;\calf_t\right],\quad \text{from}~\eqref{exmp:model},\nonumber\\
 &=& S_t\exp\left(-\sum_{k=1}^t Z_k\right)\e\left( \exp\left\{-\sum_{k=t+1}^T Z_k - T/2\right\}\left[\exp(Z_{t+1}-1/2) - 1 \right]\right)\nonumber\\
 &=& - S_0 \exp( -t) \frac{\mathrm{e}-1}{\mathrm{e}}= - 4(\mathrm{e}-1) \exp( -t-1).
\end{eqnarray} 
And obviously, since $S_t$ is a martingale under $Q_3$, we have 
\[
v_t(f_3) = \e\left[ S_{t+1} - S_t \;\vert\;\calf_t\right] = 0.
\]
Hence, for $s=(s_1,s_2,s_3)\in \real^3$, the random variable $\lambda_t(s)$ is given by
\[
\begin{split}
\lambda_t(s)&= {4\mathrm{e}^{-t}(\mathrm{e}-1)}\left[s_1\exp\left\{2\sum_{1}^t z_k \right\} - s_2\exp( - 1)\right]\\
&= {4\mathrm{e}^{-t}(\mathrm{e}-1)}\left[ s_1\mathrm{e}^t\left(\frac{S_t}{S_0}\right)^2 -  s_2\exp( - 1)\right]\\
&= 4(\mathrm{e}-1)\left[{s_1}\left(\frac{S_t}{S_0}\right)^2 - {s_2}\mathrm{e}^{-t-1}\right].
\end{split}
\]

\smallskip
Thus, for $1 \le t \le 2$ and $\mathbf{z}=(z_1,z_2,z_3)\in \Omega$, we have the following table:
\begin{align*}
v^+(f_1)(t,\mathbf{z})&=v_t(f_1)(\mathbf{z}), & d^+_1&=76.34,\qquad
& v^-(f_1)(t,\mathbf{z})&=0, & d^-_1&=0, \\
v^+(f_2)(t,\mathbf{z})&=0, & d^+_2&=0, \qquad
&v^-(f_2)(t,\mathbf{z})&=2.53\mathrm{e}^{-t}, & d^-_2&=3.80,\\
v^+(f_3)(t,\mathbf{z})&=0, & d^+_3&=0, \qquad
&v^+(f_3)(t,\mathbf{z})&=0, & d^-_3&=0.
\end{align*}
From above and \eqref{howmanyzero}, we also have $\aleph=2$.
Clearly, we need to consider only two changes of measures, the one given by $v^{+}(f_1)$ and the other by $v^-(f_2)$. The rest are all zero measures. Finally, since $a_t \equiv 0$, from \eqref{whatiscf}, we get $c(f_i) = 0, \qquad i=1,2,3$.

We take the precision parameters to be
\[
\epsilon = .5, \qquad \delta= .05.
\]
From \eqref{e:nosamples} and \eqref{e:nosamples-}, we determine a sufficient number of samples for desired accuracy would be
\[
\kappa^+_1 = 1,400,000,\qquad \kappa^-_2 = 10,500.
\]
Let us now analyze the probability measures $\mu_1^+$ and $\mu_2^{-}$ on $\rr^{3}\times\{0,1,2\}$. If $\mathbf{z} \in \rr^3$, and $0\le t \le 2$, then from \eqref{Qpm} and \eqref{exmp:whatisvtf1} we get
\begin{equation}
\begin{split}
\dx\mu_1^+(\mathbf{z},t)\; &\propto \;v^+(f_1)(\mathbf{z},t)\cdot\dx\medprob(\mathbf{z},t) \\
 & \propto \exp\left\{2\sum_1^t z_k - t\right\}\cdot \left(\frac{1}{\sqrt{2\pi}}\right)^3\exp\left\{ -\frac{1}{2}\sum_{k=1}^{3}z^2_k \right\}\\
 & \propto \mathrm{e}^{t}\left(\frac{1}{\sqrt{2\pi}}\right)^3\exp\left\{ -\frac{1}{2}\sum_{k=1}^{t}(z_k - 2)^2 -\frac{1}{2}\sum_{k=t+1}^{3}z_k^2 \right\}.
\end{split}
\end{equation}
Thus generating a sample from $\mu^+_1$ is the same as picking a $t \in (0,1,2)$ randomly with probability proportional to $\exp(t)$. Then, conditionally on $t$, we generate $t$ independent samples $Z_1,\ldots,Z_t$ from $N(2,1)$, and $3-t$ samples from $N(0,1)$.   

Simulating from $\mu^-_2$ is even simpler, since, from \eqref{exmp:whatisvtf2}, we get that
\begin{equation}
\begin{split}
\dx\mu_2^-(\mathbf{z},t)\; &\propto \;v^-(f_2)(\mathbf{z},t)\cdot\dx\medprob(\mathbf{z},t) \\
 & \propto \mathrm{e}^{- t}\cdot \left(\frac{1}{\sqrt{2\pi}}\right)^3\exp\left\{ -\frac{1}{2}\sum_{k=1}^{3}z^2_k \right\}.
\end{split}
\end{equation}
Here, we pick $t$ from $\{0,1,2\}$ with probability proportional to $\exp(-t)$, and generate $(Z_1,\ldots,Z_T)$ as independent and identically distributed samples from $N(0,1)$.

Finally, we take $\eta$ to be $N(0,1)$.
\medskip

\noindent{\sc{result of simulations.}} We first generate the required number of samples from $\mu^{+}_1$ and $\mu_2^{-}$ and set them aside. Now we choose a variety of grids, making them finer and more localized as we proceed, until $\hat{\rho}$ converges to a global minimum.
 
An estimate of the minimum capital is $w^*_0= 0.41$. The optimal grid point comes to $s_1=0.05,s_2=9.65$. Thus, an estimate of the trading strategy for this capital is $\xi^*_t = \Phi\left( \lambda_t \right)$, where $\Phi$ is the standard normal cumulative distribution function, and $\lambda_t$ is the process given by
\[
\lambda_t = 4(\mathrm{e}-1)\left[.05\left(\frac{S_t}{S_0}\right)^2 - 9.65\mathrm{e}^{-t-1}\right].
\]
In other words, with a probability more than $(1-\aleph\delta)=.9$, we will indeed have $\rho(w_0^* + W(\xi^*))\le \epsilon=0.5$.

\section{Conclusion}
We devise a Monte-Carlo algorithm to compute near-minimal initial capital requirement and a suitable trading strategy to achieve acceptability at a future date. The benefit of this approach is that it gives precise numerical values for portfolio optimization problems where purely theoretical methods (e.g. backward induction, linear programming) fail. 

The primary shortcoming is that this approach requires intensive computing, mainly due to bound \eqref{ulln}. However, the fault lies in the crudeness of the exact theoretical bound, and not in the method itself. In fact, there are better bounds (e.g. due to Talagrand \cite{talagrand}) which, unfortunately, lack exact constants.   

A related problem (brought to the authors attention by Prof. R. Jarrow at the CCCP conference, 2006) is the following. Suppose we have two risk measures $\rho_1$ and $\rho_2$. Can we find a pair $(w_0^*,\xi^*)$ of capital requirement and trading strategy,  such that $(w_0^*,\xi^*)$ minimizes $\rho_1$ among all pairs $(w,\xi)$ for which $\rho_2(w,\xi)$ is non-positive ? The author believes that the method in this paper can be suitably extended, and is currently involved in such a project.

\section{Appendix}

\subsection{Uniform law of large numbers}\label{subsec:ulln}

We briefly mention here three basic theorems about the theory of uniform law of large numbers and the related concept of Vapnik-\u{C}ervonenkis dimensions. This is a subject in itself and we shall use very little of it for our purpose. Hence we shall skip all details and refer the reader to the excellent book \cite[Chap.~12]{devroye:etal}, from where our propositions in this section have been lifted. 

\begin{notn}\label{product}We consider a probability space $\varobmu$, where $\Theta$ is a complete, separable metric space. On $\Theta^n$, let $\varrho^n$ denote the product probability measure on the product $\sigma$-algebra. Similarly on $\Theta^{\infty}:= \Theta^{\mathbb{N}}$, let $\varrho^{\infty}$ denote the infinite product probability. 
For any $\theta \in \Theta^{\infty}$, and any $n\in\mathbb{N}$, define the random \emph{empirical measure}: $\varrho_n(C):= {1}/{n} \sum_{i=1}^{n}1_{(\theta_i \in C)},\quad C \in \Im$, or, for any $\Im$-integrable function $f$, the corresponding random expectation $\varrho_n(f):= {1}/{n} \sum_{i=1}^{n}f(\theta_i)$. 
\end{notn}
\medskip

For any $C\in\Im$ and any $\epsilon > 0$, the law of large numbers dictate
\begin{subequations}
\begin{equation}
\lim_{n\rightarrow\infty}\varrho^{\infty}\Bigl( \lvert \varrho_n(C) - \varrho(C) \rvert > \epsilon\Bigr ) = 0.
\end{equation}
However, if we have a collection of $\{C_{\alpha}\}_{\alpha\in I}$ of sets in $\Im$, it is not always true that
\begin{equation}\label{uniform}
\lim_{n\rightarrow\infty}\varrho^{\infty}\Bigl( \;\sup_{\alpha \in I}\;\lvert \varrho_n(C_{\alpha}) - \varrho(C_{\alpha}) \rvert > \epsilon\Bigr ) = 0.
\end{equation}
\end{subequations}
Equality above can be achieved under proper conditions on the collection $\{C_{\alpha}\}_{\alpha \in I}$, and then we say \emph{Uniform Law of Large Numbers}(ULLN) holds. The Vapnik-\u{C}ervonenkis theory provides one such condition. Its strength lies in that the condition on $\{C_{\alpha}\}_{\alpha\in I}$ is combinatorial in nature, and hence independent from the choice of $\varrho$. (This sometimes can also be a weakness, since significant improvements can be made for specific choice of $\varrho$.) The theory begins with the concept of \emph{shatter-coefficient}. 

\begin{defn}\label{defin:shatter}
Let $\{C_{\alpha}\}_{\alpha \in I}$ be a collection of $\Im$-measurable subsets of $\Theta$. For $(\theta_1,\ldots,\theta_d)\in\Theta^d$, let $\mcal{N}(\theta_1,\ldots,\theta_d)$ be the number of different sets in 
\[
\left\{\quad\{\theta_1,\ldots,\theta_d\}\cap C_{\alpha},\; \alpha\in I\quad\right\}.
\]
The $d$-th shatter coefficient of the collection $\{C_{\alpha}\}_{\alpha \in I}$ is defined as 
\[
s_d \defin \max_{(\theta_1,\ldots,\theta_d)\in\Theta^d}\mcal{N}(\theta_1,\ldots,\theta_d).
\]
In other words, the shatter coefficient is the maximal number of different subsets of $d$ points that can be picked out by the class $\{C_{\alpha}\}_{\alpha \in I}$.
\end{defn}

\begin{rmk}
Note that we have deliberately suppressed mentioning the class $\{C_{\alpha}\}_{\alpha \in I}$ in the notation for the shatter coefficient. This is really for notational clarity. The shatter coefficient is clearly a property of the collection of sets we consider.
\end{rmk}

\noindent The following theorem can be found in \cite[~Thm 12.5, p. 197]{devroye:etal}.

\begin{thm}\label{vc}
For any collection $\{C_{\alpha}\}_{\alpha\in I}$, and for any $n\in\mathbb{N}$, $\epsilon > 0$, we have
\begin{equation}\label{eq:shatter}
\varrho^{\infty}\left\{ \sup_{\alpha\in I} \Large\lvert \varrho_n(C_{\alpha})-\varrho(C_{\alpha})\; \Large\rvert > \epsilon \right\} \leq 8s_n\exp(-n\epsilon^2/32), 
\end{equation}
where the constant $s_n$ is the $n$th shatter coefficient of the collection $\{C_{\alpha}\}_{\alpha\in I}$ and is independent of the probability measure $\varrho$. 
\end{thm}
Hence \eqref{uniform} will hold if the constants $s_n$ grows at most polynomially. This is achieved for certain collections of sets which have a finite Vapnik-\u{C}eronenkis (VC) dimension. The following definition is from \cite[~p. 196]{devroye:etal}. 

\begin{defn}\label{defn:vc}
As before we consider the collection $\{C_{\alpha}\}_{\alpha \in I}$ of $\Im$-measurable subsets of $\Theta$. The largest positive integer for which $s_d = 2^d$ is known as the VC dimension of the collection $\{C_{\alpha}\}_{\alpha \in I}$. If $s_d=2^d$ for all integers $d \ge 1$, we then define the VC dimension to be $\infty$. 
\end{defn}

The next lemma \cite[p.~218]{devroye:etal} describes a fundamental relationship between VC dimension and the shatter coefficients.
\begin{Sauers}
Let $\{C_{\alpha}\}_{\alpha\in I}$ be a subset of $\Im$ with finite VC dimension $\V > 2$. Then for all $n > 2\V$, we have $s_n \leq n^{\V}$.
\end{Sauers}
Thus Theorem~\ref{vc} together with Sauer's Leamma will yield the following.
\begin{thm}\label{vc2}
Let $(\Theta, \Im)$ be a measurable space. Let $\{C_{\alpha}\}_{\alpha\in I}$ be any collection of measurable subsets of $\Theta$ with a finite VC dimension $\V$. Then for any probability measure $\varrho$ on $(\Theta, \Im)$ and any $n\ge 2\V$, we have
\begin{equation}\label{ulln2}
\varrho^{\infty}\left\{ \sup_{\alpha\in I} \Large\lvert \varrho_n(C_{\alpha})-\varrho(C_{\alpha})\; \Large\rvert > \epsilon \right\} \leq 8n^{\V}\exp(-n\epsilon^2/32).
\end{equation}
In particular, $\lim_{n\rightarrow\infty} \varrho^{\infty}\left\{ \sup_{\alpha\in I} \Large\lvert \varrho_n(C_{\alpha})-\varrho(C_{\alpha})\; \Large\rvert > \epsilon \right\}=0$.
\end{thm}
The following better bound is from Devroye (1982).
\begin{thm}\label{combovc}
In the setting of the previous theorem \ref{vc2}, we have
\[
\varrho^{\infty}\left\{ \sup_{\alpha\in I} \Large\lvert \varrho_n(C_{\alpha})-\varrho(C_{\alpha})\; \Large\rvert > \epsilon \right\} \leq 4s_{n^2}\exp(-2n\epsilon^2 + 4\epsilon + 4\epsilon^2).
\]
Hence, by Sauer's Lemma,
\begin{equation}\label{ulln}
\varrho^{\infty}\left\{ \sup_{\alpha\in I} \Large\lvert \varrho_n(C_{\alpha})-\varrho(C_{\alpha})\; \Large\rvert > \epsilon \right\} \leq 4n^{2 \V_{C}}\exp(-2n\epsilon^2+4\epsilon + 4\epsilon^2).
\end{equation}
\end{thm}

Finally, we shall need the following collection of sets with finite VC dimension. 
\begin{prop}\label{dud}\cite[Thm 7.2]{dudley}
Let $G$ be a $d$-dimensional real vector space of real functions on an infinite set $X$. Define the class of sets 
\[
\C=\left\{ \left\{x\in X: g(x) > 0\right\}: g\in G\right\}.
\] 
Then the VC dimension of $\C$ is not more than $d$.
\end{prop}

\bibliographystyle{amsalpha}
%\bibliography{/Users/soumikpal/lib/soumiks.bib}

\providecommand{\bysame}{\leavevmode\hbox to3em{\hrulefill}\thinspace}
\providecommand{\MR}{\relax\ifhmode\unskip\space\fi MR }
\providecommand{\MRhref}[2]{%
  \href{http://www.ams.org/mathscinet-getitem?mr=#1}{#2}
}
\providecommand{\href}[2]{#2}

\end{document}